\documentclass{article}

\usepackage{amsfonts}
\usepackage{amssymb,amsmath,latexsym}
\usepackage{amsthm}
\usepackage{eurosym}
\usepackage{ragged2e}
\usepackage{float}

\newtheorem{theorem}{Theorem}

\newtheorem{corollary}{Corollary}

\newtheorem{definition}{Definition}
\newtheorem{example}{Example}

\newtheorem{lemma}{Lemma}

\newtheorem{remark}{Remark}

\numberwithin{equation}{section}

\begin{document}

\textbf{Formulas of special polynomials involving Bernoulli polynomials derived from matrix equations and Laplace transform}\\

\textbf{Ezgi Polat and Yilmaz Simsek} \\

Department of Mathematics, Faculty of Science, University of Akdeniz, Antalya TR-07058,Turkey\\

\textbf{E-mail.} ezgipolat1247@gmail.com and ysimsek@akdeniz.edu.tr\\

 \textbf{Abstract.}	The main purpose and motivation of this article is to 	create a linear transformation on the polynomial ring of rational numbers. A 	matrix representation of this linear transformation based on standard fundamentals will be given. For some special cases of this matrix, matrix 	equations including inverse matrices, the Bell polynomials will be given. 	With the help of these equations, new formulas containing different 	polynomials, especially the Bernoulli polynomials, will be given. Finally, by applying the Laplace transform to the generating function for the 	Bernoulli polynomials, we derive some novel formulas involving the Hurwitz zeta function and infinite series.\\



\section{Introduction}

It is well known that in recent years many different applications of linear
transformations have been given not only in the algebraic field of
mathematics but also in other applied sciences. To give an example of these
different applications, geometric transformations implemented in computer
graphics also occur; that is translation, rotation and scaling of $2D$ or $%
3D $ objects can be done using a transformation matrix. Therefore, linear
transformations are also used as a tool to describe change. Many examples
can be given for these, some of which are well known to be used in analysis
as transformations corresponding to derivatives, transformations corresponding to integrals, or in relativity as a device to keep track of
local transformations of reference frames. In addition, other implementation
examples can be given, including compiler optimizations of nested loop code
and parallelization of compiler techniques.

In recent years, it has also been seen that generating functions for the
Bernoulli numbers and polynomials were given by different methods. For
example, in complex analysis, generating functions of these numbers and
polynomials can be constructed using the Cauchy derivative formula, Cauchy
residue theorem and meromorphic functions. These numbers and polynomials are
also used in algebraic topology and number theory, as criteria of regular
prime numbers or their appearance in the Todd class which can be seen on
complex vector bundles in topological space, the values of zeta functions on
even integers, in Milnor's homotopy group related to the characteristic
class. They can also appear in $K$-theory, the Euler-Maclaurin summation
formula, Bessel functions, trigonometric functions, cylindrical functions,
hypergeometric functions etc. In addition, it is also well-known that these
numbers and polynomials are given by the Volkenborn integral on the set of $p
$-adic integers. Perhaps their definitions are also given on other spaces or
sets that we have not seen before. Therefore, the main motivation of this
paper is to give different computational formulas of polynomials containing
the Bernoulli polynomials using not only the linear transformation defined
on the polynomials ring of rational numbers and its matrix equations, but
also the Laplace transform and the Hurwitz zeta function.

We can use the following notations and definitions:

Let $\mathbb{N},$ $\mathbb{Q},$ $\mathbb{R}$ and $\mathbb{C}$ denote a set
of positive integers, the ring of rational integers,\ a set of real numbers,
and a set of complex numbers. $\mathbb{N}_{0}\mathbb{=N\cup }\left\{
0\right\} $.

The Bernoulli numbers and polynomials are respectively given by%
\begin{equation}
	F(t)=\frac{t}{e^{t}-1}=\sum_{n=0}^{\infty }B_{n}\frac{t^{n}}{n!}  \label{Bn}
\end{equation}%
and%
\begin{equation}
	F(t,x)=e^{tx}F(t)=\sum_{n=0}^{\infty }B_{n}(x)\frac{t^{n}}{n!}  \label{Bnx}
\end{equation}
(\textit{cf}. \cite{Arakawa,Riardon,Srivastava}).

The Stirling numbers of the first kind and the second kind, which are
denoted by $S_{1}\left( n,k\right) $ and $S_{2}(m,n)$, respectively, are
defined by 
\begin{equation}
	\frac{\left( \log \left( 1+t\right) \right) ^{k}}{k!}=\sum\limits_{n=0}^{%
		\infty }S_{1}\left( n,k\right) \frac{t^{n}}{n!}  \label{S11}
\end{equation}%
and%
\begin{equation}
	\frac{(e^{t}-1)^{n}}{n!}=\sum\limits_{m=0}^{\infty }S_{2}(m,n)\frac{t^{m}}{m!%
	} \label{sn}
\end{equation}%
(\textit{cf}. \cite{Charambides,Riardon,Srivastava}).

The array polynomials are defined by the following generating function:%
\begin{equation}
	\frac{\left( e^{t}-1\right) ^{v}e^{tx}}{v!}=\sum_{m=0}^{\infty }S_{v}^{m}(x)%
	\frac{t^{m}}{m!}  \label{sna}
\end{equation}%
(\textit{cf}. \cite{Cagic,Chang,Simsekarray}).

The results of this article, including all sections, are briefly stated as
follows:

In Section \ref{section2}, we defined $\mathbb{Q}$-linear operator on the polynomials
ring $\mathbb{Q}[x]$. By using this operator, we give its matrix
representation with respect to the basis of $\mathbb{Q}[x]$\textbf{: }$%
\left\{ 1,x,x^{2},x^{3},\ldots \right\} $. By using matrix representation,
we find some new classes of special polynomials involving Bernoulli
polynomials and Bell polynomials.

In Section \ref{section3}, we give linear transformation. By applying Cayley--Hamilton
theorem, we also give inverse matrix formulas involving the Bell polynomials
for this linear transformation. Using this inverse matrix formula, we derive
matrix representation for the Bernoulli polynomials.

In Section \ref{section4}, we give a derivative formulas for linear transformation.

In Section \ref{section5}, we define a new family of polynomials and its matrix. By using this matrix with its inverse, we also derive the
Bernoulli numbers.

In Section \ref{section6}, applying the Laplace transform to generating function for the
Bernoulli polynomials, we give not only infinite series representation for
the Bernoulli polynomials, but also we derive some novel formulas including
the Stirling numbers and the array polynomials.

Finally, we give conclusion section of this paper.

\section{A class of $\mathbb{Q}$-linear operator on the polynomials ring $%
	\mathbb{Q}[x]$ and their applications} \label{section2}

In this section, we define not only the following $\mathbb{Q}$-linear
operator on the polynomials ring $\mathbb{Q}[x]$, but also construct its
matrix representation with respect to the basis of $\mathbb{Q}[x]$\textbf{: }%
$\left\{ 1,x,x^{2},x^{3},\ldots \right\} $. By using this operator and its
matrix, we derive some new classes of special polynomials involving
the Bernoulli polynomials, the Bell polynomials etc.

Let%
\begin{equation*}
	\mathcal{E}:\mathbb{Q}[x]\rightarrow \mathbb{Q}[x].
\end{equation*}%
Let $g(x)$ be integrable function. Then, we define the following linear
operator on the polynomials ring $\mathbb{Q}[x]$:%
\begin{equation}
	\mathcal{E}\left[ g(x)\right] =\int_{ax+b}^{cx+d}g(u)du,  \label{AE-1}
\end{equation}%
where%
\begin{equation*}
	g(x)\in \mathbb{Q}[x],
\end{equation*}%
$a,b,m,j\in \mathbb{R}$.

Substituting $g(x)=x^{n}$ into (\ref{AE-1}), we obtain%
\begin{equation}
	\mathcal{E}\left[ x^{n};a,b,c,d\right] =\int_{ax+b}^{cx+d}u^{n}du=\frac{%
		\left( cx+d\right) ^{n+1}-\left( ax+b\right) ^{n+1}}{n+1}.  \label{af}
\end{equation}

\begin{remark}
	Substituting $a=c=d=1$ and $b=0$ into \textup{(\ref{af})}, $\mathcal{E}\left[
	x^{n};1,0,1,1\right] $ reduces to $\mathbb{Q}$-linear map, which is given by
	Arakawa et al. \textup{\cite[p. 55]{Arakawa}}.
\end{remark}

Using (\ref{af}), we give the following formulas for any polynomials%
\begin{equation*}
	P_{n}(x)=\sum_{j=0}^{n}\alpha _{j}x^{j}
\end{equation*}%
on the ring $\mathbb{Q}[x]$ under the linear operation $\mathcal{E}$:

\begin{eqnarray}
	&&\mathcal{E}\left[ P_{n}(x);a,b,c,d\right]  \label{AE-1a} \\
	&=&\sum_{j=0}^{n}\binom{j+1}{0}\frac{\alpha _{j}}{j+1}(d^{j+1}-b^{j+1})+%
	\sum_{j=0}^{n}\binom{j+1}{1}\frac{\alpha _{j}}{j+1}(cd^{j}-ab^{j})x \notag
	\\
	&&+\sum_{j=1}^{n}\binom{j+1}{2}\frac{\alpha _{j}}{j+1}%
	(c^{2}d^{j-1}-a^{2}b^{j-1})x^{2}+\sum_{j=2}^{n}\binom{j+1}{3}\frac{\alpha
		_{j}}{j+1}(c^{3}d^{j-2}-a^{3}b^{j-2})x^{3}  \notag \\
	&&+\sum_{j=3}^{n}\binom{j+1}{4}\frac{\alpha _{j}}{j+1}%
	(c^{4}d^{j-3}-a^{4}b^{j-3})x^{4}+\ldots \notag\\ &&+\sum_{j=k}^{n}\binom{j+1}{k-1}\frac{%
		\alpha _{j}}{j+1}(c^{k+1}d^{j-k}-a^{k+1}b^{j-k})x^{k+1}+  \notag \\
	&&\ldots +\alpha _{n-1}\frac{\binom{n}{n}}{n}(c^{n}d-a^{n}b)x^{n}+\alpha _{n}%
	\frac{\binom{n+1}{n+1}}{n+1}(c^{n+1}-a^{n+1})x^{n+1}.  \notag
\end{eqnarray}

With the aid of (\ref{AE-1a}), $(n+1)\times (n+2)$ matrix representation of
the linear operation $\mathcal{E}\left[ P_{n}(x);a,b,c,d\right] $\ with
respect to basis of $\left\{ 1,x,x^{2},x^{3},\ldots \right\} $ is
given by the following theorem:
\begin{theorem}
	\begin{eqnarray}
		&	M_{\mathcal{E}}(P_{n}(x);a,b,c,d) = \hspace{6.5cm} \label{AE-2a} \\\notag \\
		&	\begingroup\makeatletter\def\f@size{2}\check@mathfonts 	 \begin{bmatrix}
			\frac{\binom{1}{0}(d-b)}{1} & \frac{\binom{1}{1}(c-a)}{1} & 0 & 0 & \ldots & 
			0 \\ 
			\frac{\binom{2}{0}(d^{2}-b^{2})}{2} & \frac{\binom{2}{1}(cd-ab)}{2} & \frac{%
				\binom{2}{2}(c^{2}-a^{2})}{2} & 0 & \ldots & 0 \\ 
			\frac{\binom{3}{0}(d^{3}-b^{3})}{3} & \frac{\binom{3}{1}(cd^{2}-ab^{2})}{3}
			& \frac{\binom{3}{2}(c^{2}d-a^{2}b)}{3} & \frac{\binom{3}{3}(c^{3}-a^{3})}{3}
			& \ldots & 0 \\ 
			\vdots & \vdots & \vdots & \vdots & \vdots & \vdots \\ 
			\frac{\binom{k+1}{0}(d^{k+1}-b^{k+1})}{k+1} & \frac{\binom{k+1}{1}%
				(cd^{k}-ab^{k})}{k+1} & \frac{\binom{k+1}{2}(c^{2}d^{k-1}-a^{2}b^{k-1})}{k+1}
			& \frac{\binom{k+1}{3}(c^{3}d^{k-2}-a^{3}b^{k-2})}{k+1} & \ddots & 0 \\ 
			\vdots & \vdots & \vdots & \vdots & \vdots & \vdots \\ 
			\frac{\binom{n+1}{0}(d^{n+1}-b^{n+1})}{n+1} & \frac{\binom{n+1}{1}%
				(cd^{n}-ab^{n})}{n+1} & \frac{\binom{n+1}{2}(c^{2}d^{n-1}-a^{2}b^{n-1})}{n+1}
			& \frac{\binom{n+1}{3}(c^{3}d^{n-2}-a^{3}b^{n-2})}{n+1} & \ldots & \frac{%
				\binom{n+1}{n+1}(c^{n+1}-a^{n+1})}{n+1}
		\end{bmatrix} 
		\notag \endgroup
	\end{eqnarray} 	
\end{theorem}

If $c=a$, (\ref{AE-2a}) is a square matrix and $\det (M_{\mathcal{E}%
}(P_{n}(x);a,b,a,d)):=\left\vert M_{\mathcal{E}}\right\vert \neq 0$. Thus we
define the following a new family of polynomials:
\begin{definition}
	We define $Q_{n}(x)$ polynomial sequences with respect to basis of $\mathbb{Q}[x]$\textbf{: }$\left\{ 1,x,x^{2},x^{3},\ldots \right\} $:
	
	\begin{equation}
		\begin{bmatrix}
			Q_{0}(x) \\ 
			Q_{1}(x) \\ 
			Q_{2}(x) \\ 
			\vdots  \\ 
			Q_{n-1}(x)%
		\end{bmatrix}%
		=\left( \left( M_{\mathcal{E}}(P_{n}(x);a,b,a,d)\right) ^{-1}\right) 
		\begin{bmatrix}
			1 \\ 
			x \\ 
			x^{2} \\ 
			\vdots  \\ 
			x^{n-1}%
		\end{bmatrix}%
		, \label{Me}
	\end{equation}%
	where $A^{-1}$ denotes inverse of the matrix $A$.
\end{definition}

\begin{example}
	\label{Example 1}Putting $n=2$ in \textup{(\ref{Me})}, we get%
	\begin{equation*}
		M_{\mathcal{E}}(P_{2}(x);a,b,c,d)=%
		\begin{bmatrix}
			d-b & c-a & 0 \\ 
			\frac{d^{2}-b^{2}}{2} & cd-ab & c^{2}-a^{2}%
		\end{bmatrix}%
		\text{,}
	\end{equation*}%
	\begin{eqnarray*}
		\left\vert M_{\mathcal{E}}(P_{2}(x);a,b,a,d)\right\vert &=&a\left(
		b-d\right) ^{2} \\
		&=&ab^{2}-2abd+ad^{2}
	\end{eqnarray*}%
	and%
	\begin{equation*}
		\left( M_{\mathcal{E}}(P_{2}(x);a,b,a,d)\right) ^{-1}=%
		\begin{bmatrix}
			\frac{-1}{b-d} & 0 \\ 
			\frac{\left( b+d\right) }{2(ab-ad)} & \frac{-1}{ab-ad}%
		\end{bmatrix}%
		\text{.}
	\end{equation*}%
	We have%
	\begin{equation*}
		\begin{bmatrix}
			Q_{0}(x) \\ 
			Q_{1}(x)%
		\end{bmatrix}%
		=%
		\begin{bmatrix}
			\frac{-1}{b-d} & 0 \\ 
			\frac{\left( b+d\right) }{2a(b-d)} & \frac{-1}{a(b-d)}%
		\end{bmatrix}%
		\begin{bmatrix}
			1 \\ 
			x%
		\end{bmatrix}%
		\text{.}
	\end{equation*}%
	From the above equation, we get the following polynomials:%
	\begin{eqnarray*}
		Q_{0}(x) &=&-\frac{1}{b-d} \\
		Q_{1}(x) &=&\frac{\left( b+d\right) }{2a(b-d)}-\frac{1}{a(b-d)}x\text{.}
	\end{eqnarray*}
\end{example}

\begin{example}
	\label{Example 2}Putting $n=3$ in \textup{(\ref{Me})}, we get%
	\begin{equation*}
		M_{\mathcal{E}}(P_{3}(x);a,b,c,d)=%
		\begin{bmatrix}
			(d-b) & (c-a) & 0 & 0 \\ 
			\frac{(d^{2}-b^{2})}{2} & (cd-ab) & \frac{(c^{2}-a^{2})}{2} & 0 \\ 
			\frac{(d^{3}-b^{3})}{3} & (cd^{2}-ab^{2}) & (c^{2}d-a^{2}b) & \frac{%
				(c^{3}-a^{3})}{3}%
		\end{bmatrix}%
		.
	\end{equation*}
\end{example}

By using inverse matrix method we get%
\begin{equation*}
	\left( M_{\mathcal{E}}(P_{3}(x);a,b,a,d)\right) ^{-1}=%
	\begin{bmatrix}
		\frac{-1}{b-d} & 0 & 0 \\ 
		\frac{(b+d)}{2a(b-d)} & \frac{-1}{a(b-d)} & 0 \\ 
		\frac{-\left( \frac{b^{2}+d^{2}}{6}+\frac{2bd}{3}\right) }{a^{2}(b-d)} & 
		\frac{\left( b+d\right) }{a^{2}(b-d)} & \frac{-1}{a^{2}(b-d)}%
	\end{bmatrix}.
\end{equation*}%
We have%
\begin{equation*}
	\begin{bmatrix}
		Q_{0}(x) \\ 
		Q_{1}(x) \\ 
		Q_{2}(x)%
	\end{bmatrix}%
	=%
	\begin{bmatrix}
		\frac{-1}{b-d} & 0 & 0 \\ 
		\frac{(b+d)}{2a(b-d)} & \frac{-1}{a(b-d)} & 0 \\ 
		\frac{-\left( \frac{b^{2}+d^{2}}{6}+\frac{2bd}{3}\right) }{a^{2}(b-d)} & 
		\frac{\left( b+d\right) }{a^{2}(b-d)} & \frac{-1}{a^{2}(b-d)}%
	\end{bmatrix}%
	\begin{bmatrix}
		1 \\ 
		x \\ 
		x^{2}%
	\end{bmatrix}%
	.
\end{equation*}%
Therefore from the above, we get the following polynomials:%
\begin{eqnarray*}
	Q_{0}(x) &=&\frac{-1}{b-d} \\
	Q_{1}(x) &=&\frac{(b+d)}{2a(b-d)}-\frac{1}{a(b-d)}x \\
	Q_{2}(x) &=&\frac{-\left( \frac{b^{2}+d^{2}}{6}+\frac{2bd}{3}\right) }{%
		a^{2}(b-d)}+\frac{\left( b+d\right) }{a^{2}(b-d)}x+\frac{-1}{a^{2}(b-d)}%
	x^{2}.
\end{eqnarray*}

\section{Applications of the linear transformation $\mathcal{E}\left[
	P_{n}(x);a,b,c,d\right] $}\label{section3}

In this section, we give some special values of the linear transformation $%
\mathcal{E}\left[ P_{n}(x);a,b,c,d\right] $. By applying Cayley--Hamilton
theorem, we also give inverse matrix formulas involving the Bell
polynomials. By the aid of the inverse matrix formula, we give matrix
representation for the Bernoulli polynomials.

Substituting $a=c$ into (\ref{AE-1a}), we get%
\begin{eqnarray*}
	&&\mathcal{E}\left[ P_{n}(x);a,b,a,d\right] \\
	&=&\sum_{j=0}^{n}\binom{j+1}{0}\frac{\alpha _{j}}{j+1}(d^{j+1}-b^{j+1})+a%
	\sum_{j=0}^{n}\binom{j+1}{1}\frac{\alpha _{j}}{j+1}(d^{j}-b^{j})x+ \\
	&&+a^{2}\sum_{j=1}^{n}\binom{j+1}{2}\frac{\alpha _{j}}{j+1}%
	(d^{j-1}-b^{j-1})x^{2}+a^{3}\sum_{j=2}^{n}\binom{j+1}{3}\frac{\alpha _{j}}{%
		j+1}(d^{j-2}-b^{j-2})x^{3} \\
	&&+\ldots +a^{n}\alpha _{n-1}\frac{\binom{n+1}{n}}{n+1}(d-b)x^{n}.
\end{eqnarray*}

Substituting $a=c$ and $d=b+1$ into (\ref{AE-1a}), we also get%
\begin{eqnarray}
	&&\mathcal{E}\left[ P_{n}(x);a,b,a,b+1\right]   \label{AEf} \\
	&=&\sum_{j=0}^{n}\binom{j+1}{0}\frac{\alpha _{j}}{j+1}\left( \left(
	(b+1)^{j+1}-b^{j+1}\right) +a(j+1)((b+1)^{j}-b^{j})x\right)   \notag \\
	&&+a^{2}\sum_{j=1}^{n}\binom{j+1}{2}\frac{\alpha _{j}}{j+1}\left(
	(b+1)^{j-1}-b^{j-1}\right) x^{2}+\ldots   \notag \\
	&&+a^{n}\sum_{j=n-1}^{n}\binom{j+1}{n}\frac{\alpha _{j}}{j+1}\left(
	(b+1)^{j-n+1}-b^{j-n+1}\right) x^{n}
	.
\end{eqnarray}%
We assume that $j-n+1\geq 0$, otherwise when $j-n+1<0$, these power's
numbers omitted in the related sums.

With the aid of (\ref{AEf}), $n+1\times n+1$ matrix representation of the
linear operation $\mathcal{E}\left[ P_{n}(x);a,b,a,b+1\right] $\ with
respect to basis of \textbf{\ }$\left\{ 1,x,x^{2},x^{3},\ldots \right\} $ is
given by the following corollary:

\begin{corollary}
	\begin{eqnarray}
		M_{\mathcal{E}}(P_{n}(x);a,b,a,b+1)= \hspace{8.5cm}\label{AC-MT} \\\notag \\
		\begingroup\makeatletter\def\f@size{9}\check@mathfonts
		\begin{bmatrix}
			1 & 0 & 0 & \ldots & 0 \\ 
			\frac{\binom{2}{0}}{2}((b+1)^{2}-b^{2}) & a & 0 & \ldots & 0 \\ 
			\frac{\binom{3}{0}}{3}((b+1)^{3}-b^{3}) & \frac{\binom{3}{1}}{3}%
			a((b+1)^{2}-b^{2}) & a^{2} & \ldots & 0 \\ 
			\vdots & \vdots & \vdots & \ddots & \ldots \\ 
			\frac{\binom{n+1}{0}}{n+1}((b+1)^{n+1}-b^{n+1}) & \frac{\binom{n+1}{1}}{n}%
			a((b+1)^{n}-b^{n}) & \frac{\binom{n+1}{2}}{n}a^{2}((b+1)^{n-1}-b^{n-1}) & 
			\ldots & a^{n}%
		\end{bmatrix}.
		\notag  \endgroup
	\end{eqnarray}
\end{corollary}

By applying Cayley--Hamilton theorem to (\ref{AC-MT}) and using inverse
matrix formula, given in formula (\textit{cf.} \cite[Eq. (B.10) and Eq. (B.11)]%
{Kondratyuk tersmatris}), we show that inverse matrix of the matrix $M_{%
	\mathcal{E}}(P_{n}(x);a,b,a,d)$ is given by the following theorem.

\begin{lemma}
	If $c=a$, $M_{\mathcal{E}}(P_{n}(x);a,b,c,d)$ is a diagonal matrix, then we
	have%
	\begin{eqnarray}
		&&\left( M_{\mathcal{E}}(P_{n}(x);a,b,a,d)\right) ^{-1}  \label{im-1} \\
		&=&\frac{1}{\left\vert M_{\mathcal{E}}\right\vert }\sum_{j=0}^{n}M_{\mathcal{%
				E}}^{j}(P_{n}(x);a,b,a,d)\notag\\
		&&\times	\sum_{a_{1},a_{2},\ldots ,a_{n}}\prod_{v=1}^{n}%
		\frac{(-1)^{a_{v}+1}}{v^{a_{v}}a_{v}!}\left( tr(M_{\mathcal{E}%
		}^{v}(P_{n}(x);a,b,a,d)\right) ^{a_{v}}  \notag
	\end{eqnarray}%
	$n+1$ is a dimension of matrix $M_{\mathcal{E}}(P_{n}(x);a,b,a,d)$, $%
	tr\left\{ M_{\mathcal{E}}(P_{n}(x);a,b,a,d)\right\} $ is the trace of main
	diagonal matrix $M_{\mathcal{E}}(P_{n}(x);a,b,a,d)$, and also%
	\begin{equation*}
		\sum_{a_{1},a_{2},a_{3},\ldots ,a_{n}}
	\end{equation*}%
	means that this sum is taken over $j$ and the sets of all $a_{v}\geq 0$
	satisfying the following equation;%
	\begin{equation*}
		j+\sum_{m=0}^{n}m.a_{m}=n.
	\end{equation*}
\end{lemma}

\begin{lemma}
	If $c=a$, $M_{\mathcal{E}}(P_{n}(x);a,b,c,d)$ is a diagonal matrix, then we
	have%
	\begin{eqnarray}
		&&\left( M_{\mathcal{E}}(P_{n}(x);a,b,a,d)\right) ^{-1}  \label{im-2} \\
		&=&\frac{1}{\left\vert M_{\mathcal{E}}\right\vert }\sum_{j=1}^{n}M_{\mathcal{%
				E}}^{j-1}(P_{n}(x);a,b,a,d)\frac{(-1)^{n-1}}{(n-j)!}B_{n-j}(w_{1},w_{2},%
		\ldots ,w_{n-j}),  \notag
	\end{eqnarray}%
	where%
	\begin{equation*}
		w_{k}=-(k-1)!tr\left\{ M_{\mathcal{E}}^{k}(P_{n}(x);a,b,a,d)\right\} 
	\end{equation*}%
	and $B_{n}(w_{1},w_{2},\ldots ,w_{n})$ denotes the $n$th complete
	exponential Bell polynomials, defined by \textup{(see \cite{Bell,Charambides})}
	\begin{equation*}
		B_{n}(w_{1},w_{2},\ldots ,w_{n})=n!\sum_{n=1k_{1}+2k_{2}+\ldots
			+nk_{n}}\prod\limits_{j=1}^{n}\frac{w_{j}^{k_{j}}}{\left( j!\right)
			^{k_{j}}k_{j}!}.
	\end{equation*}
\end{lemma}

By using (\ref{im-1}) and $\mathcal{E}^{-1}\left[ x^{n};a,b,a,b+1\right] $,
we find another family of polynomials by the following theorem:

\begin{theorem}
	\begin{equation}
		\begin{bmatrix}
			H_{0}(x) \\ 
			H_{1}(x) \\ 
			H_{2}(x) \\ 
			\vdots \\ 
			H_{n-1}(x)%
		\end{bmatrix}%
		=\left( \left( M_{\mathcal{E}}(P_{n}(x);a,b,a,b+1)\right) ^{-1}\right) 
		\begin{bmatrix}
			1 \\ 
			x \\ 
			x^{2} \\ 
			\vdots \\ 
			x^{n-1}%
		\end{bmatrix}.
		\label{AC-MT2}
	\end{equation}
\end{theorem}

When $a=c=d=1$ and $b=0$, the Bernoulli polynomials, which are also given by
the following corollary:

\begin{corollary}
	\begin{equation}
		\begin{bmatrix}
			B_{0}(x) \\ 
			B_{1}(x) \\ 
			B_{2}(x) \\ 
			\vdots \\ 
			B_{n-1}(x)%
		\end{bmatrix}%
		=\left( \left( M_{\mathcal{E}}(P_{n}(x);1,0,1,1)\right) ^{-1}\right) 
		\begin{bmatrix}
			1 \\ 
			x \\ 
			x^{2} \\ 
			\vdots \\ 
			x^{n-1}%
		\end{bmatrix}%
		,  \label{AC-MT1}
	\end{equation}%
	where 
	\begin{eqnarray*}
		&&\left( M_{\mathcal{E}}(P_{n}(x);1,0,1,1)\right) ^{-1} \\
		&=&\sum_{j=0}^{n}\frac{M_{\mathcal{E}}^{j}(P_{n}(x);1,0,1,1)}{\det (M_{%
				\mathcal{E}}(P_{n}(x);1,0,1,1))}\\
		&&\times \sum_{a_{1},a_{2},\ldots
			,a_{n}}\prod_{v=1}^{n}\frac{(-1)^{a_{v}+1}}{v^{a_{v}}a_{v}!}\left( tr(M_{%
			\mathcal{E}}^{v}(P_{n}(x);1,0,1,1)^{a_{v}}\right) .
	\end{eqnarray*}
\end{corollary}

\begin{remark}
	When $a=c=d=1$ and $b=0$, with the aid of \textup{(\ref{AC-MT1})}, we see that 
	\begin{equation*}
		B_{n}(x)=\mathcal{E}^{-1}\left[ x^{n};1,0,1,1\right]
	\end{equation*}%
	(cf. \textup{\cite[p. 55]{Arakawa}}).
\end{remark}

\section{Derivative formula for linear transformation $\mathcal{E}\left[
	P_{n}(x);a,b,c,d\right] $}\label{section4}

In this section, we give derivative formula of the linear transformation $%
\mathcal{E}\left[ P_{n}(x);a,b,c,d\right] $

Taking derivative of $\mathcal{E}\left[ P_{n}(x);a,b,c,d\right] $ with
respect to $x$, we arrive at the following derivative formula:%
\begin{eqnarray*}
	&&\frac{d}{dx}\left\{ \mathcal{E}\left[ P_{n}(x);a,b,c,d\right] \right\} \\
	&&=\sum_{j=0}^{n}\binom{j+1}{1}\frac{\alpha _{j}(cd^{j}-b^{j}a)}{j+1}%
	+2\sum_{j=1}^{n}\binom{j+1}{2}\frac{(c^{2}d^{j-1}-a^{2}b^{j-1})\alpha _{j}}{%
		j+1}x \\
	&&+3\sum_{j=2}^{n}\binom{j+1}{3}\frac{(c^{3}d^{j-2}-a^{3}b^{j-2})\alpha _{j}%
	}{j+1}x^{2}+4\sum_{j=3}^{n}\binom{j+1}{4}\frac{(c^{4}d^{j-3}-a^{4}b^{j-3})%
		\alpha _{j}}{j+1}x^{3} \\
	&&+\ldots +\sum_{j=k}^{n}\binom{j+1}{k-1}\frac{\left( k+1\right)
		(c^{k+1}d^{j-k}-a^{k+1}b^{j-k})\alpha _{j}}{j+1}x^{k} \\
	&&+\ldots +(c^{n}d-a^{n}b)\alpha _{n-1}x^{n-1}+(c^{n+1}-a^{n+1})\alpha
	_{n}x^{n}.
\end{eqnarray*}

Putting $a=c=d=1$ and $b=0$ in\ the above equation, we have%
\begin{eqnarray*}
	&&\frac{d}{dx}\left\{ \mathcal{E}\left[ P_{n}(x);1,0,1,1\right] \right\}  \\
	&=&\sum_{j=0}^{n}\binom{j+1}{1}\frac{\alpha _{j}}{j+1}+2\sum_{j=1}^{n}\binom{%
		j+1}{2}\frac{\alpha _{j}}{j+1}x+3\sum_{j=2}^{n}\binom{j+1}{3}\frac{\alpha
		_{j}}{j+1}x^{2} \\
	&&+4\sum_{j=3}^{n}\binom{j+1}{4}\frac{\alpha _{j}}{j+1}x^{3}+\ldots
	+\sum_{j=k}^{n}\binom{j+1}{k-1}\frac{\left( k+1\right) \alpha _{j}}{j+1}%
	x^{k}+\ldots +\alpha _{n-1}x^{n-1}.
\end{eqnarray*}
Combining (\ref{AC-MT1}) with the above equation and its matrix, we have the
following well-known result for the Bernoulli polynomials:%
\begin{equation*}
	\frac{d}{dx}\left\{ B_{n}(x)\right\} =nB_{n-1}(x).
\end{equation*}

\section{A new family of polynomials and its matrix}\label{section5}

In this section, we define a new family of polynomials. By using these
polynomials, we also define a matrix representation of the coefficients of
these polynomials. We show that special values of these matrix with its
inverse produce the Bernoulli numbers.

By using (\ref{af}), we define the following polynomial:%
\begin{eqnarray*}
	Y_{n}(x;a,b,c,d) &=&\left( cx+d\right) ^{n+1}-\left( ax+b\right) ^{n+1} \\
	&=&\sum\limits_{j=0}^{n+1}\binom{n+1}{j}(c^{j}d^{n+1-j}-a^{j}b^{n+1-j})x^{j}.
\end{eqnarray*}%
Some values of $Y_{n}(x;a,b,c,d)$ are given by%
\begin{eqnarray*}
	Y_{0}(x;a,b,c,d) &=&(c-a)x+d-b \\
	Y_{1}(x;a,b,c,d) &=&(c^{2}-a^{2})x^{2}+2(cd-ab)x+(d^{2}-b^{2}) \\
	Y_{2}(x;a,b,c,d)
	&=&(c^{3}-a^{3})x^{3}+3(c^{2}d-a^{2}b)x^{2}+3(cd^{2}-ab^{2})x+(d^{3}-b^{3})
	\\
	&&\vdots  \\
	Y_{n}(x;a,b,c,d) &=&(c^{n+1}-a^{n+1})x^{n+1}+\binom{n+1}{n}%
	(c^{n}d-a^{n}b)x^{n}+\cdots \\&&+(d^{n+1}-b^{n+1}).
\end{eqnarray*}%
We define matrix representation of the aid of coefficients of the above
polynomials as follows:
\begin{eqnarray}
	M(Y_{n}(x;a,b,c,d))= \hspace{8.5cm}\label{Y_m}
	\\\notag \\
	\begingroup\makeatletter\def\f@size{5}\check@mathfonts
	\begin{bmatrix}
		d-b & c-a & 0 & 0 & \ldots & 0 \\ 
		d^{2}-b^{2} & 2(cd-ab) & c^{2}-a^{2} & 0 & \ldots & 0 \\ 
		d^{3}-b^{3} & 3(cd^{2}-ab^{2}) & 3(c^{2}d-a^{2}b) & c^{3}-a^{3} & \ldots & 0
		\\ 
		\vdots & \vdots & \vdots & \vdots & \vdots & \vdots \\ 
		d^{n+1}-b^{n+1} & \binom{n+1}{1}\left( cd^{n}-ab^{n}\right) & \binom{n+1}{2}%
		\left( c^{2}d^{n-1}-a^{2}b^{n-1}\right) & \binom{n+1}{3}\left(
		c^{3}d^{n-2}-a^{3}b^{n-2}\right) & \ldots & c^{n+1}-a^{n+1}%
	\end{bmatrix}.
	\endgroup
	\notag
\end{eqnarray}
Thus we get%
\begin{equation*}
	\left[ 
	\begin{array}{c}
		Y_{0}(x;a,b,a,d) \\ 
		Y_{1}(x;a,b,a,d) \\ 
		Y_{2}(x;a,b,a,d) \\ 
		\vdots \\ 
		Y_{n-1}(x;a,b,a,d)%
	\end{array}%
	\right] =M^{-1}(Y_{n}(x;a,b,a,d))%
	\begin{bmatrix}
		1 \\ 
		x \\ 
		x^{2} \\ 
		\vdots \\ 
		x^{n-1}%
	\end{bmatrix}.
\end{equation*}

The first row of the inverse of the matrix $M(Y_{n}(x;a,b,a,d))$ given by
equation (\ref{Y_m}) can also be considered as the generating function for a
special family of numbers, including Bernoulli numbers maybe other certain
family of special numbers. That is, substituting $a=c=d=1$ and $b=0$ into
inverse matrix $M^{-1}(Y_{n}(x;1,0,1,1))$, given by equation (\ref{Y_m}),
all entries of the matrix $M^{-1}(Y_{n}(x;1,0,1,1))$ are reduced to the
Bernoulli numbers, respectively.

Some well-known examples are given as follows:%
\begin{equation*}
	M(Y_{6}(x;1,0,1,1))=%
	\begin{bmatrix}
		1 & 0 & 0 & 0 & 0 & 0 & 0 \\ 
		1 & 2 & 0 & 0 & 0 & 0 & 0 \\ 
		1 & 3 & 3 & 0 & 0 & 0 & 0 \\ 
		1 & 4 & 6 & 4 & 0 & 0 & 0 \\ 
		1 & 5 & 10 & 10 & 5 & 0 & 0 \\ 
		1 & 6 & 15 & 20 & 15 & 6 & 0 \\ 
		1 & 7 & 21 & 35 & 35 & 21 & 7%
	\end{bmatrix}.
\end{equation*}%
By using (\ref{im-1}) inverse matrix method we get%
\begin{equation*}
	M^{-1}(Y_{6}(x;1,0,1,1))=%
	\begin{bmatrix}
		1 & 0 & 0 & 0 & 0 & 0 & 0 \\ 
		-\frac{1}{2} & \frac{1}{2} & 0 & 0 & 0 & 0 & 0 \\ 
		\frac{1}{6} & -\frac{1}{2} & \frac{1}{3} & 0 & 0 & 0 & 0 \\ 
		0 & \frac{1}{4} & -\frac{1}{2} & \frac{1}{4} & 0 & 0 & 0 \\ 
		-\frac{1}{30} & 0 & \frac{1}{3} & -\frac{1}{2} & \frac{1}{5} & 0 & 0 \\ 
		0 & -\frac{1}{12} & 0 & \frac{5}{12} & -\frac{1}{2} & \frac{1}{6} & 0 \\ 
		\frac{1}{42} & 0 & -\frac{1}{6} & 0 & \frac{1}{2} & -\frac{1}{2} & \frac{1}{7%
		}%
	\end{bmatrix}%
	,
\end{equation*}%
where we obtained Bernoulli numbers up to $n=6$ in the first column of the $%
M(Y_{6}(x;a,b,a,d))$ matrix (cf.\ for detail, A027642, see also \cite{Helms}%
). On the other hand, other columns or rows of the above matrix also
represent different number families. These number sequences may also be in
the class of known other families of special number. For this, it is
recommended to examine the \textit{oeis} page (https://oeis.org/).

\section{Applications of the Laplace transform to generating function for
	the Bernoulli polynomials}\label{section6}

In this section, applications of the Laplace transform and the Mellin
transformation to (\ref{Bnx}), we give a relation between the Bernoulli
polynomials and the Hurwitz zeta function. We also give infinite series
representation involving the Bernoulli polynomials.

Let $x\in 
\mathbb{C}
$ with $x=a+ib$\ and $\overline{x}=a-ib$ with $a>0$ and $b>0$. We set%
\begin{equation*}
	F(-u,\overline{x})=\sum_{n=0}^{\infty }(-1)^{n}B_{n}(a-ib)\frac{u^{n}}{n!}.
\end{equation*}%
Therefore%
\begin{equation}
	\frac{ue^{-ua}}{1-e^{-u}}=e^{-uib}\sum_{n=0}^{\infty }(-1)^{n}B_{n}(a-ib)%
	\frac{u^{n}}{n!}.  \label{bx11}
\end{equation}%
Eq. (\ref{bx11}) can also be written as%
\begin{equation}
	e^{-ua}\sum_{n=1}^{\infty }\frac{\left( 1-e^{-u}\right) ^{n-1}}{n}%
	=e^{-uib}\sum_{n=0}^{\infty }(-1)^{n}B_{n}(a-ib)\frac{u^{n}}{n!}.  \label{Ey}
\end{equation}%
Combining the above equation with (\ref{sn}), we get%
\begin{eqnarray*}
	&&\sum_{m=0}^{\infty }\sum\limits_{v=0}^{m}\binom{m}{v}\sum_{n=1}^{v}\frac{%
		(-1)^{m+n-1}(n-1)!a^{m-v}S_{2}(v,n-1)}{n}\frac{u^{m}}{m!} \\
	&=&\sum_{m=0}^{\infty }(-1)^{m}\sum\limits_{v=0}^{m}\binom{m}{v}%
	(ib)^{m-v}B_{v}(a-ib)\frac{u^{m}}{m!}.
\end{eqnarray*}%
Substituting the following well-known formula into the right-hand side of
the above equation%
\begin{equation*}
	\sum\limits_{v=0}^{m}\binom{m}{v}x^{m-v}B_{v}(y)=B_{m}(x+y),
\end{equation*}%
we obtain%
\begin{eqnarray*}
	&&\sum_{m=0}^{\infty }\sum\limits_{v=0}^{m}\binom{m}{v}\sum_{n=1}^{v}\frac{%
		(-1)^{m+n-1}(n-1)!a^{m-v}S_{2}(v,n-1)}{n}\frac{u^{m}}{m!} \\
	&=&\sum_{m=0}^{\infty }(-1)^{m}B_{m}(a)\frac{u^{m}}{m!}.
\end{eqnarray*}%
By comparing the coefficients $\frac{u^{m}}{m!}$\ on both sides of the above
equation, we arrive at the following theorem:

\begin{theorem}
	Let $m\in \mathbb{N}$. Then we have%
	\begin{equation}
		B_{m}(a)=\sum\limits_{v=1}^{m}\binom{m}{v}\sum_{n=1}^{v}\frac{%
			(-1)^{n-1}(n-1)!a^{m-v}S_{2}(v,n-1)}{n}.  \label{Ey1}
	\end{equation}
\end{theorem}

By combining (\ref{sna}) with (\ref{Ey}), by the same method of Cakic and
Milovanovic \cite{Cagic}, Chang and Ha \cite{Chang}, Simsek\cite{Simsekarray}%
, we have%
\begin{equation*}
	\sum_{m=0}^{\infty }\sum_{n=1}^{\infty }\frac{\left( -1\right) ^{n-1}(n-1)!}{%
		n}S_{n-1}^{m}(a)\frac{\left( -u\right) ^{m}}{m!}=\sum_{m=0}^{\infty
	}\sum\limits_{v=0}^{m}\binom{m}{v}(ib)^{m-v}B_{v}(a-ib)\frac{\left(
		-u\right) ^{m}}{m!}.
\end{equation*}%
After some calculation, we have the following result, which also proved in 
\cite{Chang}:%
\begin{equation}
	\sum_{n=1}^{m}\frac{\left( -1\right) ^{n-1}\left( n-1\right) !}{n}%
	S_{n-1}^{m}(x)=\sum\limits_{v=0}^{m}\binom{m}{v}(ib)^{m-v}B_{v}(a-ib).
	\label{af-I2}
\end{equation}

Combining the above equation with (\ref{Ey1}) yields%
\begin{equation}
	\sum\limits_{v=1}^{m}\sum_{n=1}^{v}\frac{(-1)^{n-1}\binom{m}{v}%
		(n-1)!a^{m-v}S_{2}(v,n-1)}{n}=\sum_{n=0}^{m}\frac{\left( -1\right) ^{n}n!}{%
		n+1}S_{n}^{m}(x).  \label{af-I0}
\end{equation}%
Combining (\ref{af-I2}) and (\ref{af-I0}) with%
\begin{equation}
	\frac{\left( -1\right) ^{n}n!}{n+1}=\sum\limits_{v=0}^{n}S_{1}(n,v)B_{v},
	\label{dn2}
\end{equation}%
(\textit{cf}. \cite{Charambides,KimDah,Riardon,SimsekRev}),
we arrive at the following theorem:
\begin{theorem}
	Let $m\in \mathbb{N}_{0}$. Then we have%
	\begin{equation}
		B_{m}(a)=\sum_{n=0}^{m}\sum\limits_{v=0}^{n}S_{1}(n,v)B_{v}S_{n}^{m}(a)
		\label{af-I01}
	\end{equation}%
	and%
	\begin{equation*}
		\sum\limits_{v=0}^{m}\binom{m}{v}(ib)^{m-v}B_{v}(a-ib)=\sum_{n=0}^{m}\sum%
		\limits_{v=0}^{n}S_{1}(n,v)B_{v}S_{n}^{m}(a).
	\end{equation*}
\end{theorem}

Here, we note that different proof of the equation (\ref{af-I0}) was also
given by Chang and Ha \cite{Chang}.

By applying the Laplace transform to generating functions involving
Bernoulli polynomials, recently many interesting studies have been published
(\textit{cf}. \cite{Revista gun simsek,goss,Laplace,Srivastava}). Therefore, by using (\ref{bx11}), we have%
\begin{equation*}
	\sum_{n=0}^{\infty }ue^{-u(a+n)}=e^{-uib}\sum_{n=0}^{\infty
	}(-1)^{n}B_{n}(a-ib)\frac{u^{n}}{n!}\text{.}
\end{equation*}%
Integrate equation the above equation with respect to $u$ from $1$ to $%
\infty $, we get%
\begin{equation}
	\sum_{n=0}^{\infty }\int_{0}^{\infty }ue^{-u(a+n)}du=\sum_{n=0}^{\infty }%
	\frac{(-1)^{n}}{n!}B_{n}(a-ib)\int_{0}^{\infty }u^{n}e^{-ibu}du  \label{bx0}
\end{equation}%
with $a>0$. Applying\ the Laplace transform of the function $g(u)=u^{n}$:%
\begin{equation*}
	\mathcal{L}\left\{ g(u)\right\} =\frac{n!}{y^{n+1}}
\end{equation*}%
(where $y>0$) to the both sides of (\ref{bx0}), we obtain the following
theorem:

\begin{theorem}
	Let $b>0$ and $a>0$. Then we have%
	\begin{equation}
		\zeta (2,a)=\sum_{n=0}^{\infty }\frac{i^{n-1}B_{n}(a-ib)}{b^{n+1}}.
		\label{bx13}
	\end{equation}
\end{theorem}

By using (\ref{bx11}), we have%
\begin{equation*}
	ue^{-au}=\left( e^{-ibu}-e^{-(1+ib)u}\right) \sum_{n=0}^{\infty
	}(-1)^{n}B_{n}(a-ib)\frac{u^{n}}{n!}.
\end{equation*}%
Applying\ the Laplace transform to the above equation, for $a>0$ and $b>0$,
we get%
\begin{eqnarray*}
	\int_{0}^{\infty }ue^{-au}du &=&\sum_{n=0}^{\infty }(-1)^{n}B_{n}(a-ib)\frac{%
		1}{n!}\int_{0}^{\infty }u^{n}e^{-ibu}du \\
	&&-\sum_{n=0}^{\infty }(-1)^{n}B_{n}(a-ib)\frac{1}{n!}\int_{0}^{\infty
	}u^{n}e^{-(1+ib)u}du.
\end{eqnarray*}%
After some calculations, for $b>1$, we arrive at the following theorem:

\begin{theorem}
	Let $a>0$ and $b>1$. Then we have%
	\begin{equation}
		\sum_{n=0}^{\infty }\sum_{j=0}^{n}(-1)^{n}\binom{n}{j}\frac{\left( \left(
			1+ib\right) ^{n+1}-\left( ib\right) ^{n+1}\right) (a-ib)^{n-j}B_{j}}{\left(
			ib-b^{2}\right) ^{n+1}}=\frac{1}{a^{2}}.  \label{af-I0a}
	\end{equation}
\end{theorem}

\begin{remark}
	Putting $n=1$ into equation \textup{(17)} in \textup{\cite{Revista gun simsek}}, a
	relationship can be established with equation \textup{(\ref{af-I0a})}.
\end{remark}

Substituting $a=1$ into (\ref{bx13}), and combining with the following known
formula%
\begin{equation*}
	\zeta (2):=\zeta (2,1)=\sum_{n=1}^{\infty }\frac{1}{n^{2}}=\frac{\pi ^{2}}{6}%
	,
\end{equation*}%
we get the following corollary:

\begin{corollary}
	Let $b>1$ and $a>0$. Then we have%
	\begin{equation*}
		\sum_{n=0}^{\infty }\frac{i^{n-1}B_{n}(a-ib)}{b^{n+1}}=\frac{\pi ^{2}}{6}.
	\end{equation*}
\end{corollary}

Replacing $t,\ x$ by $-t$, $y-x$ in (\ref{Bnx}) yields%
\begin{equation*}
	F(-t,y-x)=\frac{-te^{-t(y-x)}}{e^{-t}-1}=\sum_{n=0}^{\infty
	}(-1)^{n}B_{n}(y-x)\frac{t^{n}}{n!}.
\end{equation*}%
Taking the partial derivative $k$ times with respect to $y$ in the above
equation, we get%
\begin{equation*}
	\frac{\partial ^{k}}{\partial y^{k}}\left\{ F(-t,y-x)\right\} =\frac{%
		(-1)^{k+1}t^{k+1}e^{-t(y-x)}}{e^{-t}-1}.
\end{equation*}%
After some calculations, we obtain%
\begin{equation*}
	\frac{(-1)^{k+1}t^{k+1}e^{-ty}}{e^{-t}-1}=e^{-tx}\sum_{n=0}^{\infty }\frac{%
		(-1)^{n}\frac{\partial ^{k}}{\partial y^{k}}\left\{ B_{n}(y-x)\right\} }{n!}%
	t^{n}.
\end{equation*}%
Integrate equation the above equation with respect to $u$ from $1$ to $%
\infty $, we get%
\begin{equation*}
	(-1)^{k+1}\int_{0}^{\infty }\frac{t^{k+1}e^{-ty}}{e^{-t}-1}%
	dt=\sum_{n=0}^{\infty }\frac{(-1)^{n}\frac{\partial ^{k}}{\partial y^{k}}%
		\left\{ B_{n}(y-x)\right\} }{n!}\int_{0}^{\infty }t^{n}e^{-tx}dt,
\end{equation*}%
which yields%
\begin{equation*}
	(-1)^{k+1}(k+1)!\sum_{n=0}^{\infty }\frac{1}{(y+n)^{k+2}}=\sum_{n=0}^{\infty
	}\frac{(-1)^{n}}{x^{n+1}}\frac{\partial ^{k}}{\partial y^{k}}\left\{
	B_{n}(y-x)\right\} .
\end{equation*}%
After some calculations, we arrive at the following theorem:

\begin{theorem}
	Let $k\in \mathbb{N}_{0}$. Let $y>0$. Then we have%
	\begin{equation}
		\zeta (k+2,y)=\frac{(-1)^{k+1}}{(k+1)!}\sum_{n=0}^{\infty }\frac{(-1)^{n}}{%
			x^{n+1}}\frac{\partial ^{k}}{\partial y^{k}}\left\{ B_{n}(y-x)\right\} .
		\label{bnx4}
	\end{equation}
\end{theorem}

\begin{description}
	\item Substituting $k=0$ into (\ref{bnx4}), we arrive at (\ref{bx13}).
\end{description}

Let the ratio of two convergent series be given such that when the first
term of the series in the denominator is different from zero, using (\ref{Bn}%
), we have the following result:%
\begin{equation*}
	\frac{1+0t+0t^{2}+\cdots }{1+\frac{1}{2!}t+\frac{1}{3!}t^{2}+\frac{1}{4!}%
		t^{3}+\cdots }=\sum_{n=0}^{\infty }B_{n}\frac{t^{n}}{n!}.
\end{equation*}%
From the above equation, we get%
\begin{equation*}
	1=\left( 1+\frac{1}{2!}t+\frac{1}{3!}t^{2}+\frac{1}{4!}t^{3}+\cdots \right)
	\left( B_{0}+\frac{B_{1}}{1!}t+\frac{B_{2}}{2!}t^{2}+\frac{B_{3}}{3!}%
	t^{3}+\cdots \right) .
\end{equation*}%
Using the Cauchy product rule on the right-hand side of the above equation
yields an infinite system of linear equations with unknown $B_{k}$. So, this
system is of special form since every $k\in \mathbb{N}_{0}$ the first $k+1$
equations contain only the first $k+1$ unknown $B_{k}$. Therefore, the
solution of the system of linear equations is given by the following
determinant, which is calculated the Bernoulli numbers outside the influence
of their generating function:%
\begin{equation*}
	B_{n}=(-1)^{n}n!%
	\begin{vmatrix}
		\frac{1}{2!} & 1 & 0 & 0 & 0 & \ldots  & 0 \\ 
		\frac{1}{3!} & \frac{1}{2!} & 1 & 0 & 0 & \ldots  & 0 \\ 
		\frac{1}{4!} & \frac{1}{3!} & \frac{1}{2!} & 1 & 0 & \ldots  & 0 \\ 
		\vdots  & \vdots  & \vdots  & \vdots  & \vdots  & \ddots  & 0 \\ 
		\frac{1}{(n+1)!} & \frac{1}{n!} & \frac{1}{(n-1)!} & \frac{1}{(n-2)!} & 
		\frac{1}{(n-3)!} & \ldots  & \frac{1}{2!}%
	\end{vmatrix}%
\end{equation*}
(\textit{cf.} \cite[p. 149]{E.Pap}, \cite{H.Chen}).

\section{Conclusion}

In this paper, we defined $\mathbb{Q}$-linear operator on the polynomials
ring $\mathbb{Q}[x]$. By using this operator, we gave its matrix
representation with respect to the basis of $\mathbb{Q}[x]$\textbf{: }$%
\left\{ 1,x,x^{2},x^{3},\ldots \right\} $. Using matrix representation,
we obtanied some new classes of special polynomials involving Bernoulli
polynomials and Bell polynomials. We also defined other linear transformations. By applying Cayley--Hamilton
theorem, we gave their inverse matrix formulas involving the Bell polynomials. Using this inverse matrix formula, we gave
matrix representation for the Bernoulli polynomials. Moreover, applying the Laplace transform to generating function for the
Bernoulli polynomials, we obtanied both infinite series representation for
the Bernoulli polynomials and many new formulas associated with
the Stirling numbers and the array polynomials. 

The results of this article may potentially  be used both in applied sciences and in many branches of mathematics.

\end{document}